\documentclass{ifacconf}

\usepackage{graphicx}
\usepackage{natbib}
\usepackage{amsmath}
\usepackage{amssymb}
\usepackage{mathtools}
\usepackage{enumitem}
\usepackage{tikz}
\usepackage{booktabs}
\usepackage{array}

\setlist[enumerate]{label={\textnormal{(\alph*)}}, ref={(\alph*)}, leftmargin=*, nosep}

\DeclareMathOperator{\Real}{Re} 
\DeclareMathOperator{\Imag}{Im} 
\DeclareMathOperator{\id}{Id} 
\DeclareMathOperator{\diff}{d\!} 

\DeclarePairedDelimiter{\abs}{\lvert}{\rvert} 
\DeclarePairedDelimiter{\norm}{\lVert}{\rVert} 

\newcommand{\suchthat}{\ifnum\currentgrouptype=16 \mathrel{}\middle|\mathrel{}\else\mid\fi}

\begin{document}

\begin{frontmatter}

\title{On qualitative properties of single-delay linear retarded differential equations: Characteristic roots of maximal multiplicity are necessarily dominant\thanksref{footnoteinfo}} 

\thanks[footnoteinfo]{Cor\-res\-pon\-ding author: Guilherme Mazanti (gui\-lher\-me.ma\-zan\-ti@l2s.cen\-tra\-le\-sup\-elec.fr)}

\author[LSS,IPSA]{Guilherme Mazanti} 
\author[LSS,IPSA]{Islam Boussaada} 
\author[LSS]{Silviu-Iulian Niculescu}

\address[LSS]{University Paris-Saclay, CNRS, CentraleSup\'elec, Laboratoire des Signaux et Syst\`emes (L2S), Inria Saclay, DISCO Team, 3, rue Joliot Curie, 91192, Gif-sur-Yvette, France.}
\address[IPSA]{Institut Polytechnique des Sciences Avanc\'ees (IPSA)\\63 boulevard de Brandebourg, 94200 Ivry-sur-Seine, France.}
\address{\{first-name.last-name\}@l2s.centralesupelec.fr}

\begin{abstract}
This paper presents necessary and sufficient conditions for the existence of a real root of maximal multiplicity in the spectrum of a linear time-invariant single-delay equation of retarded type. We also prove that this root is always strictly dominant, and hence determines the asymptotic behavior of the system. These results are based on improved a priori bounds on the imaginary part of roots on the complex right half-plane.
\end{abstract}


\end{frontmatter}

\section{Introduction}

\setlength{\parskip}{6pt plus 2pt minus 1pt} 

Systems with time delays are useful models in a wide range of scientific and technological domains such as biology, chemistry, economics, physics, or engineering, the presence of the delays usually coming from propagation phenomena, such as of material, energy, or information, with a finite propagation speed. Due to their numerous applications, these kinds of systems have been the subject of much attention by researchers in several fields, in particular since the 1950s and 1960s. We refer to \cite{Gu2003Stability, Hale1993Introduction, Michiels2014Stability} for details on time-delay systems and their applications.

In this paper, we consider linear time-invariant equations with a single delay of the form
\begin{equation}
\label{MainSystTime}
\begin{split}
y^{(n)}(t) & + a_{n-1} y^{(n-1)}(t) + \dotsb + a_0 y(t) \\
 & + \alpha_{n-1} y^{(n-1)}(t - \tau) + \dotsb + \alpha_0 y(t - \tau) = 0,
\end{split}
\end{equation}
where $n$ is a positive integer, the coefficients $a_k$ and $\alpha_k$ are real numbers for $k \in \{0, \dotsc, n-1\}$, and the delay $\tau$ is a positive real number. Equation \eqref{MainSystTime} is said to be of \emph{retarded type} since the derivative of highest order only appears in the non-delayed term $y^{(n)}(t)$.

The stability analysis of equations of the form \eqref{MainSystTime} and more general time-delay systems has attracted much research effort and is an active field (see, e.g., \cite{Chen1995New, Gu2003Stability, Michiels2014Stability, Olgac2002Exact, Sipahi2011Stability, Li2017Frequency}). In the delay-free situation, i.e., when $\alpha_k = 0$ for every $k \in \{0, \dotsc, n-1\}$, the stability of \eqref{MainSystTime} can be studied through spectral methods by considering the corresponding characteristic polynomial, whose complex roots determine the asymptotic behavior of solutions of the system. Under the presence of delays, the asymptotic behavior of \eqref{MainSystTime} is also determined by the complex roots of some characteristic function $\Delta: \mathbb C \to \mathbb C$ (see, e.g., \cite{Hale1993Introduction, Michiels2014Stability, Mori1982Estimate}), defined for $s \in \mathbb C$ by
\begin{equation}
\label{Delta}
\Delta(s) = s^n + \sum_{k=0}^{n-1} a_k s^k + e^{-s\tau} \sum_{k=0}^{n-1} \alpha_k s^k.
\end{equation}
Functions of the form \eqref{Delta} are particular instances of quasipolynomials, defined as follows (see, e.g., \cite{Wielonsky2001Rolle, Berenstein1995Complex}).

\begin{defn}
A \emph{quasipolynomial} is an entire function $Q$ which can be written under the form
\[
Q\left(s; \lambda_0, (a_{j0})_{j=0}^{d_0}, \dotsc, \lambda_\ell, (a_{j\ell})_{j=0}^{d_\ell}\right) = \sum_{k=0}^\ell p_k(s) e^{-\lambda_k s},
\]
where $\ell$ is a positive integer, $\lambda_0, \dotsc, \lambda_\ell$ are pairwise distinct real numbers, and, for $k \in \{0, \dotsc, \ell\}$, $a_{d_k k} \neq 0$ and $p_k$ is the polynomial $p_k(s) = \sum_{j=0}^{d_k} a_{jk} s^j$ of degree $d_k$. The integer $D = \ell + \sum_{k=0}^\ell d_k$ is called the \emph{degree} of $Q$.
\end{defn}

When $\lambda_0 = 0$ and $\lambda_1, \dotsc, \lambda_\ell$ are positive, the above general quasipolynomial is the characteristic equation of a linear time-invariant delayed equation with delays $\lambda_1, \dotsc, \lambda_\ell$. A classical result on quasipolynomials provided in \cite[Problem 206.2]{Polya1998Problems} (which holds independently on the sign of $\lambda_k$), known as the \emph{P\'{o}lya--Szeg\H{o} bound}, implies that, given a quasipolynomial $Q$ of degree $D \geq 0$, the multiplicity of any root of $Q$ does not exceed $D$. For the quasipolynomial $\Delta$ from \eqref{Delta}, this means that any of its roots has multiplicity at most $2n$. Recent works such as \cite{Boussaada2016Characterizing, Boussaada2016Tracking} have provided characterizations of multiple roots of quasipolynomials using approaches based on Birkhoff and Vandermonde matrices.

Similarly to the delay-free case, the quasipolynomial $\Delta$ is useful for the stability analysis of \eqref{MainSystTime}. Indeed, all solutions of \eqref{MainSystTime} converge exponentially fast to $0$ if and only if $\Real s < 0$ for every $s \in \mathbb C$ such that $\Delta(s) = 0$, and the asymptotic behavior of solutions of \eqref{MainSystTime} is determined by the real number $\gamma_0 = \sup\{\Real s \suchthat s \in \mathbb C,\; \Delta(s) = 0\}$, called the \emph{spectral abscissa} of $\Delta$. However, contrarily to the delay-free case, \eqref{Delta} has infinitely many roots. Designing a general criterion in terms of the coefficients $a_k, \alpha_k$ for $k \in \{0, \dotsc, n-1\}$ and the delay $\tau$ ensuring that all roots of $\Delta$ have negative real part is a difficult problem (see \cite{Hayes1950Roots} for some early results in the case $n = 1$), which is in contrast with the delay-free case in which Routh--Hurwitz criterion provides a stability test in terms of the coefficients of the characteristic polynomial.

The spectral abscissa of $\Delta$ is related to the notion of dominant root, defined as follows.
\begin{defn}
Let $Q: \mathbb C \to \mathbb C$ and $s_0 \in \mathbb C$ be such that $Q(s_0) = 0$. We say that $s_0$ is a \emph{dominant} (respectively, \emph{strictly dominant}) root of $Q$ if, for every $s \in \mathbb C \setminus \{s_0\}$ such that $Q(s) = 0$, one has $\Real s \leq \Real s_0$ (respectively, $\Real s < \Real s_0$).
\end{defn}
Dominant roots may not exist for a given function of a complex variable, but they always exist for functions of the form \eqref{Delta}, as a consequence, for instance, of the fact that $\Delta$ has finitely many roots on any vertical strip in the complex plane (see, e.g., \cite[Chapter 1, Lemma 4.1]{Hale1993Introduction}). Notice also that exponential stability of \eqref{MainSystTime} is equivalent to the dominant root of $\Delta$ having negative real part.

It turns out that, for quasipolynomials, the notions of real roots of high multiplicity and dominance are often related (see, e.g., \cite{Boussaada2016Tracking, BoussaadaMultiplicity, Boussaada2018Further}), in the sense that real roots of high multiplicity tend to be dominant, a property usually referred to as \emph{multiplicity-induced dominance} (MID for short). MID has been shown to hold, for instance, in the case $n = 2$ and $\alpha_1 = 0$ in \cite{Boussaada2018Further}, proving dominance thanks to a suitable factorization of $\Delta$ when it admits a root of multiplicity $3$; and in the case $n = 2$ and $\alpha_1 \neq 0$ in \cite{BoussaadaMultiplicity}, using Cauchy's argument principle to prove dominance of the multiple root.

Another motivation for considering roots of high multiplicity is the fact that, for delay-free systems with an affine constraint on their coefficients, the spectral abscissa is minimized on a polynomial with a single root of maximal multiplicity (see \cite{Blondel2012Explicit, Chen1979Output}). Similar properties of spectral abscissa minimization on multiple roots have also been observed for some time-delay systems in \cite{Michiels2002Continuous, Ramirez2016Design, Vanbiervliet2008Nonsmooth}. Hence, the interest in investigating multiple roots does not rely on the multiplicity itself, but rather on its connection with the dominance of this root, and the corresponding applications in stability analysis and control design.

The main goal of this paper is to investigate whether MID holds for the quasipolynomial $\Delta$ from \eqref{Delta} when assigning a root with maximal multiplicity. More precisely, the questions we address in this paper are the following.
\begin{enumerate}[label={\textbf{(Q\arabic*)}}, ref={(Q\arabic*)}, nosep, leftmargin=*]
\item\label{Ques1} Is it possible to choose $a_0, \dotsc, a_{n-1}, \alpha_0, \dotsc, \alpha_{n-1} \in \mathbb R$ in such a way that a given real number $s_0$ is a root of maximal multiplicity $2n$ of $\Delta$?
\item\label{Ques2} Under the above choice, is $s_0$ (strictly) dominant?
\end{enumerate}
Our main result, Theorem~\ref{MainTheoN}, provides affirmative answers to both questions. Question~\ref{Ques1} can be addressed in a straightforward manner, whereas the answer to \ref{Ques2} relies on a sharp a priori bound on the imaginary part of roots of $\Delta$ with nonnegative real part, established in Lemma~\ref{OnlyLemma}.

The paper is organized as follows. Notations used in the paper are standard. Section~\ref{SecN2} provides a detailed study of \eqref{MainSystTime} and \eqref{Delta} in the case $n = 2$, its main result being necessary and sufficient conditions for the existence of a root of maximal multiplicity and the proof of strict dominance of this root. Notice that this case has already been considered in \cite{BoussaadaMultiplicity}, but using a different strategy for addressing \ref{Ques2}, based on an application of Cauchy's argument principle to a suitable contour. As it will be clear below, our approach allows for a simpler and easier dominance proof. Our main result is presented in Section~\ref{SecMain}, which extends the analysis of Section~\ref{SecN2} to any positive integer $n$. For the brevity of the paper, the full proof is not provided here, but it can be obtained by the same arguments as those provided in Section~\ref{SecN2} and will be presented in details in an extended version of this paper. Finally, Section~\ref{SecExpl} presents an example illustrating our main result.

\section{Motivating example: Second-order case}
\label{SecN2}

Consider \eqref{Delta} with $n = 2$, i.e.,
\begin{equation}
\label{DeltaN2}
\Delta(s) = s^2 + a_1 s + a_0 +  e^{-s\tau} (\alpha_1 s + \alpha_0),
\end{equation}
and recall that the degree of $\Delta$ is at most $4$. The main result we will prove concerning this quasipolynomial is the following.

\begin{thm}
\label{MainTheo2}
Consider the quasipolynomial $\Delta$ given by \eqref{DeltaN2} and let $s_0 \in \mathbb R$.
\begin{enumerate}
\item\label{ItemA2} The number $s_0$ is a root of multiplicity $4$ of $\Delta$ if and only if the coefficients $a_0, a_1, \alpha_0, \alpha_1$, the root $s_0$ and the delay $\tau$ satisfy the relations
\begin{subequations}
\label{Coeffs2}
\begin{align}
a_1 & = -\frac{4}{\tau} - 2 s_0, & \qquad a_0 & = \frac{6}{\tau^2} + \frac{4}{\tau} s_0 + s_0^2, \displaybreak[0] \label{Coeffs2A} \\
\alpha_1 & = -\frac{2}{\tau} e^{s_0 \tau}, & \alpha_0 & = \frac{2}{\tau} e^{s_0 \tau} \left(s_0 - \frac{3}{\tau}\right). \label{Coeffs2Alpha}
\end{align}
\end{subequations}
\item\label{ItemB2} If \eqref{Coeffs2} is satisfied, then $s_0$ is a strictly dominant root of $\Delta$.
\end{enumerate}
\end{thm}

\begin{rem}
The expressions of $a_1, a_0, \alpha_1, \alpha_0$ in \eqref{Coeffs2} are singular with respect to $\tau$ as $\tau \to 0$. If one is interested in studying the behavior of the roots of $\Delta$ as $\tau \to 0$ when \eqref{Coeffs2} is satisfied, one may consider instead the quasipolynomial $s \mapsto \tau^2 \Delta(s)$, which has the same roots as $\Delta$ but whose coefficients are regular with respect to $\tau$.
\end{rem}

Notice that, up to a translation and a scaling of the spectrum represented by the change of variables $z = \tau(s - s_0)$, one may reduce to the case $s_0 = 0$ and $\tau = 1$, in which \eqref{Coeffs2} reduces to
\[
a_1 = -4, \quad a_0 = 6, \quad \alpha_1 = -2, \quad \alpha_0 = -6,
\]
yielding the quasipolynomial
\begin{equation}
\label{DeltaHat}
\widehat\Delta(z) = z^2 - 4z + 6 - e^{-z} (2 z + 6).
\end{equation}
An important ingredient for the proof of \eqref{MainTheo2} is the following result, which provides an a priori bound on the imaginary part of the roots of $\widehat\Delta$ on the closed right half-plane.

\begin{lem}
\label{OnlyLemma}
Let $z_0$ be a root of $\widehat\Delta$ with $\Real z_0 \geq 0$. Then $\abs{\Imag z_0} < 2 \pi$.
\end{lem}

\begin{pf}
Let $\sigma = \Real z_0$, $\omega = \Imag z_0$, and
\begin{equation}
\label{MatA0A1}
A_0 = \begin{pmatrix}
0 & 1 \\
-6 & 4 \\
\end{pmatrix}, \qquad A_1 = \begin{pmatrix}
0 & 0 \\
6 & 2 \\
\end{pmatrix}.
\end{equation}
Notice that $\widehat\Delta(z) = \det(z\id - A_0 - A_1 e^{-z})$ for every $z \in \mathbb C$. Hence the root $z_0$ of $\widehat\Delta$ is an eigenvalue of the matrix $A_0 + A_1 e^{-z_0}$, and thus $\abs{z_0} \leq \rho(A_0 + A_1 e^{-z_0})$. Since, by Gelfand's formula, one has $\rho(M) = \inf_{n \in \mathbb N} \norm{M^n}^{1/n}$ for every square matrix $M$ and every submultiplicative norm $\norm{\cdot}$, one obtains in particular that
\begin{equation}
\label{BoundNormFro}
\omega^4 \leq \norm{(A_0 + A_1 e^{-z_0})^2}_2^2,
\end{equation}
where $\norm{\cdot}_2$ denotes the Frobenius norm. Then
\[
\begin{split}
\omega^4 \leq {} & 36 \abs{e^{-z_0} - 1}^2 + 4 \abs{e^{-z_0} + 2}^2 \\
& {} + 144 \abs{e^{-z_0} - 1}^2 \abs{e^{-z_0} + 2}^2 \\
& {} + \abs{(2 e^{-z_0} + 4)^2 + 6 e^{-z_0} - 6}^2,
\end{split}
\]
from where we obtain that
\begin{equation}
\label{MainBound}
\begin{split}
0 \leq {} & - 992 e^{- 2 \sigma} \cos^{2}{\omega} + (464 e^{- 2 \sigma} - 192) e^{-\sigma} \cos\omega \\
& {} + 728 + 1164 e^{- 2 \sigma} + 160 e^{- 4 \sigma} - \omega^4.
\end{split}
\end{equation}
The right-hand side of the above inequality can be seen as a real polynomial of degree $2$ in the variable $e^{-\sigma} \cos\omega$ with negative leading coefficient, and hence the above inequality is satisfied only if the discriminant of the corresponding second-degree polynomial is non-negative, i.e.,
\[
(464 e^{- 2 \sigma} - 192)^2 + 3968(728 + 1164 e^{- 2 \sigma} + 160 e^{- 4 \sigma} - \omega^4) \geq 0,
\]
which implies that
\[
\omega^4 \leq 728 + 1164 e^{- 2 \sigma} + 160 e^{- 4 \sigma} + \frac{(464 e^{- 2 \sigma} - 192)^2}{3968},
\]
and, since $\sigma \geq 0$, we obtain the bound
\[
\abs{\omega} \leq \sqrt[4]{\frac{64190}{31}} < 6.75.
\]
Assume, to obtain a contradiction, that $2\pi \leq \abs{\omega} < 6.75 $. Then $\cos\omega \in (0.893, 1]$ and, using further that $\abs{464 e^{-2\sigma} - 192} \leq 272$ and $e^{-\sigma} \leq 1$, it follows from \eqref{MainBound} that $\abs{\omega} \leq \sqrt[4]{1532.94} < 2 \pi$, yielding the required contradiction. Hence $\abs{\omega} < 2 \pi$.
\end{pf}

\begin{rem}
The inequality
\begin{equation}
\label{IneqRho}
\abs{z_0} \leq \rho(A_0 + A_1 e^{-z_0})
\end{equation}
on roots $z_0$ of $\widehat\Delta$ used in the proof of Lemma~\ref{OnlyLemma} is a classical inequality for time-delay systems of retarded type (see, e.g., \cite[Proposition~1.10]{Michiels2014Stability}). However, manipulating and extracting information directly from this inequality is usually a difficult task since it involves the computation of the spectral radius of $A_0 + A_1 e^{-z_0}$ as a function of $z_0$. Several works use instead the fact that the spectral radius is a lower bound on any submultiplicative matrix norm $\norm{\cdot}$ and consider, instead of \eqref{IneqRho}, the coarser inequality
\begin{equation}
\label{IneqNorm}
\abs{z_0} \leq \norm{A_0 + A_1 e^{-z_0}},
\end{equation}
choosing typically easier-to-compute norms, such as the $1$ norm, the Frobenius norm, or the $\infty$ norm. However, \eqref{IneqNorm} is usually much coarser than \eqref{IneqRho}. As an illustration, Table~\ref{TabNorms} provide numerical values for the best bounds on $\abs{\Imag z_0}$ that one may have using \eqref{IneqRho} and \eqref{IneqNorm} for $z_0 \in \mathbb C$ with $\Real z_0 \geq 0$.

\begin{table}[ht]
\centering
\caption{Numerical bounds on $\abs{\Imag z_0}$ for complex numbers $z_0$ with $\Real z_0 \geq 0$ satisfying \eqref{IneqRho} and \eqref{IneqNorm} for three different choices of norms.}
\label{TabNorms}
\begin{tabular}{cc}
\toprule
Inequality & Numerical bound on $\abs{\Imag z_0}$ \\
\midrule
\eqref{IneqRho} & 5.9763 \\
\eqref{IneqNorm} with the $1$ norm & 10.4520 \\
\eqref{IneqNorm} with the Frobenius norm & 10.6304 \\
\eqref{IneqNorm} with the $\infty$ norm & 11.4720 \\
\bottomrule
\end{tabular}
\end{table}

It follows from the analysis of Table~\ref{TabNorms} that it is not possible to obtain the conclusion of Lemma~\ref{OnlyLemma} only from \eqref{IneqNorm} with the above choices of norms. This is what motivates the strategy used in the proof of Lemma~\ref{OnlyLemma}, which relies on Gelfand's formula for the spectral radius to obtain from \eqref{IneqRho} the inequalities
\begin{equation}
\label{IneqNormN}
\abs{z_0}^n \leq \norm{(A_0 + A_1 e^{-z_0})^n},
\end{equation}
for any submultiplicative norm $\norm{\cdot}$ and any $n \in \mathbb N$. For $n = 1$, \eqref{IneqNormN} reduces to \eqref{IneqNorm} but, for larger $n$, \eqref{IneqNormN}, while still coarser than \eqref{IneqRho}, is usually sharper than \eqref{IneqNorm} (even though this may not always be the case; see, e.g., \cite[Chapter 1, Remark 4.2]{Kato1995Perturbation}). We represent in Figure~\ref{FigNorms2} the boundary of the sets of $z_0 \in \mathbb C$ satisfying inequality \eqref{IneqRho} involving the spectral radius (solid line) and inequality \eqref{IneqNormN} with $n = 2$ for the $1$, Frobenius, and $\infty$ norms (dashed, dotted, and dash-dotted lines), using the matrices $A_0$ and $A_1$ from \eqref{MatA0A1}. The sets themselves are on the left-hand side of the represented boundaries. We also represent, in Table~\ref{TabNorms2}, numerical bounds on $\abs{\Imag z_0}$ obtained using \eqref{IneqNormN} with $n = 2$ and same choices of norms as before.

\begin{figure*}[ht]
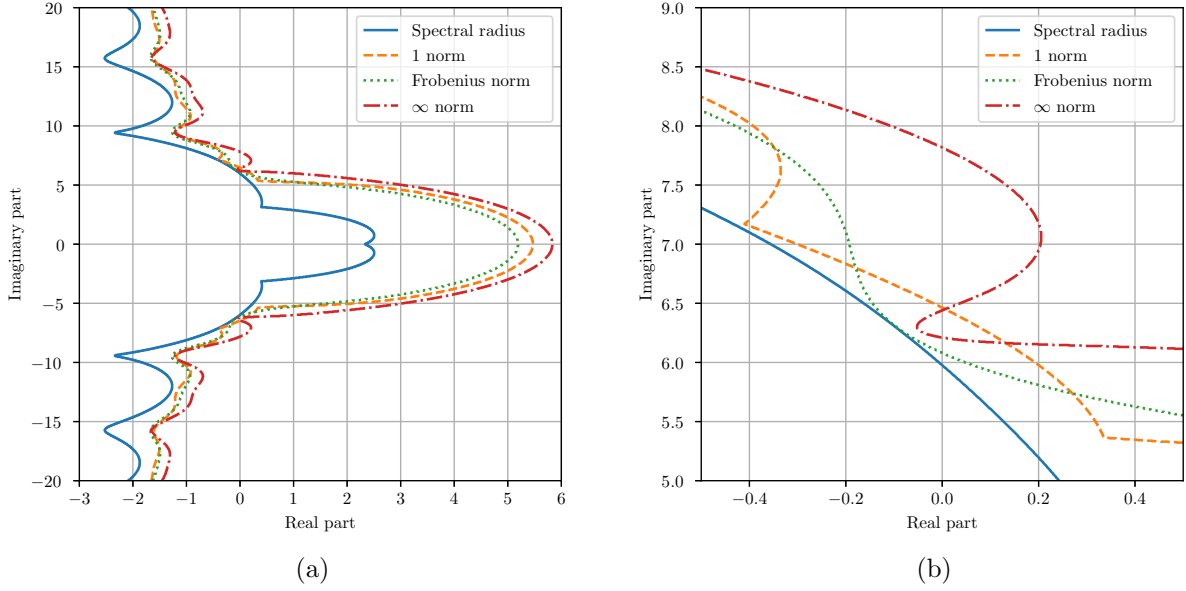

\centering
\begin{tabular}{@{} c @{} c @{}}
\resizebox{0.4494\textwidth}{!}{\input{Figures/bounds2.pgf}} & \resizebox{0.4494\textwidth}{!}{\input{Figures/bounds3.pgf}} \tabularnewline
(a) & (b)
\end{tabular}
\caption{(a) Boundaries of the sets of $z_0 \in \mathbb C$ satisfying \eqref{IneqRho} and \eqref{IneqNormN} for $n = 2$ and three different choices of norms. (b) Detail of the curves from (a) in the region $\{z \in \mathbb C \suchthat \abs{\Real z} \leq 0.5 \text{ and } 5 \leq \Imag z \leq 9\}$.}
\label{FigNorms2}
\end{figure*}

\begin{table}[ht]
\centering
\caption{Numerical bounds on $\abs{\Imag z_0}$ for complex numbers $z_0$ with $\Real z_0 \geq 0$ satisfying \eqref{IneqNormN} for $n = 2$ and three different choices of norms.}
\label{TabNorms2}
\begin{tabular}{cc}
\toprule
Inequality & Numerical bound on $\abs{\Imag z_0}$ \\
\midrule
\eqref{IneqNormN} with the $1$ norm & 6.4630 \\
\eqref{IneqNormN} with the Frobenius norm & 6.0803 \\
\eqref{IneqNormN} with the $\infty$ norm & 7.8163 \\
\bottomrule
\end{tabular}
\end{table}

Comparing Table~\ref{TabNorms} with Table~\ref{TabNorms2}, one verifies that the bounds obtained from \eqref{IneqNormN} with $n = 2$ are sharper than the ones from \eqref{IneqNorm}, and \eqref{IneqNormN} with $n = 2$ and the Frobenius norm is sufficient to prove Lemma~\ref{OnlyLemma}, justifying the strategy of our proof.

In addition to the above strategies based on \eqref{IneqRho} and inequalities between the spectral radius and matrix norms, there are also other techniques to obtain a priori bounds on roots of quasipolynomials of the form \eqref{Delta}, such as those used in \cite{Mori1989Stability} and their improved versions in \cite{Wang1992Further} and \cite{Tissir1996Further}. As a consequence of the bound by T.~Mori and H.~Kokame, one obtains that, if $z_0 \in \mathbb C$ is a root of $\widehat\Delta$ with nonnegative real part, then
\begin{equation}
\label{BoundMoriKokame}
\abs{\Imag z_0} \leq \mu(-i A_0) + \norm{A_1},
\end{equation}
where $\norm{\cdot}$ denotes a matrix norm induced by a vector norm and $\mu$ is defined for a square matrix $M$ by $\mu(M) = \lim_{\varepsilon \to 0} \frac{1}{\varepsilon}(\norm{\id + \varepsilon M} - 1)$. The improved bound by E.~Tissir and A.~Hmamed implies that
\begin{equation}
\label{BoundTissirHmamed}
\abs{\Imag z_0} \leq \mu(-i A_0) + \max_{0 \leq \theta \leq 2 \pi} \mu(A_1 e^{i\theta}).
\end{equation}
Numerical values of these bounds with $A_0$ and $A_1$ given by \eqref{MatA0A1} and using the $1$, $2$, and $\infty$ norms are presented in Table~\ref{TabOtherBounds} (note that one can no longer use the Frobenius norm since it is not induced by any vector norm). We remark that none of these bounds can be used to obtain the conclusion of Lemma~\ref{OnlyLemma}.

\begin{table}[ht]
\centering
\caption{Bounds on $\abs{\Imag z_0}$ obtained from \eqref{BoundMoriKokame} and \eqref{BoundTissirHmamed} with three different choices of norms}
\label{TabOtherBounds}
\begin{tabular}{cc}
\toprule
Inequality & Bound on $\abs{\Imag z_0}$ \\
\midrule
\eqref{BoundMoriKokame} with the $1$ norm & 12 \\
\eqref{BoundMoriKokame} with the $2$ norm & 9.8246 \\
\eqref{BoundMoriKokame} with the $\infty$ norm & 14 \\
\eqref{BoundTissirHmamed} with the $1$ norm & 12 \\
\eqref{BoundTissirHmamed} with the $2$ norm & 7.6623 \\
\eqref{BoundTissirHmamed} with the $\infty$ norm & 14 \\
\bottomrule
\end{tabular}
\end{table}

The above a priori bounds are far from forming an exhaustive list. Obtaining sharper bounds is an active research topic with important practical implications (see, e.g., \cite{CardeliquioStability} for recent developments on this topic).
\end{rem}

\begin{pf*}{Proof of Theorem~\ref{MainTheo2}.}
Let us denote by $\widetilde\Delta$ the quasipolynomial obtained by multiplying $\Delta$ by $\tau^2$ and performing the change of variables $z = \tau (s - s_0)$, i.e.,
\begin{equation}
\label{RelationDeltaTildeDelta}
\widetilde\Delta(z) = \tau^2 \Delta\left(s_0 + \tfrac{z}{\tau}\right) = z^2 + b_1 z + b_0 + e^{-z} (\beta_1 z + \beta_0),
\end{equation}
where
\begin{equation}
\label{RelationsBA}
\begin{aligned}
b_1 & = (a_1 + 2 s_0) \tau, & b_0 & = (s_0^2 + a_1 s_0 + a_0) \tau^2, \\
\beta_1 & = \alpha_1 \tau e^{- s_0 \tau}, & \beta_0 & = (\alpha_0 + \alpha_1 s_0) \tau^2 e^{- s_0 \tau}.
\end{aligned}
\end{equation}
Note that $s_0$ is a root of multiplicity $4$ of $\Delta$ if and only if $0$ is a root of multiplicity $4$ of $\widetilde\Delta$. Recalling that the degree of the quasipolynomial $\widetilde\Delta$ is $4$, $0$ is a root of multiplicity $4$ if and only if $\widetilde\Delta^{(k)}(0) = 0$ for $k \in \{0, 1, 2, 3\}$. We compute
\begin{align*}
\widetilde\Delta^\prime(z) & = 2 z + b_1 + e^{-z} (- \beta_1 z - \beta_0 + \beta_1), \displaybreak[0] \\
\widetilde\Delta^{\prime\prime}(z) & = 2 + e^{-z} (\beta_1 z + \beta_0 - 2 \beta_1), \displaybreak[0] \\
\widetilde\Delta^{\prime\prime\prime}(z) & = e^{-z} (-\beta_1 z - \beta_0 + 3 \beta_1),
\end{align*}
and one then obtains that $0$ is a root of multiplicity $4$ of $\widetilde\Delta$ if and only if
\begin{align*}
b_0 + \beta_0 & = 0, & b_1 - \beta_0 + \beta_1 & = 0, \\
\beta_0 - 2 \beta_1 & = -2, & - \beta_0 + 3 \beta_1 & = 0.
\end{align*}
The unique solution of the above system is given by $b_1 = -4$, $b_0 = 6$, $\beta_1 = -2$, and $\beta_0 = -6$. Using \eqref{RelationsBA}, one immediately verifies that these conditions are equivalent to \eqref{Coeffs2}, concluding the proof of \ref{ItemA2}. Moreover, under these conditions, one has $\widetilde\Delta = \widehat\Delta$, where $\widehat\Delta$ is defined in \eqref{DeltaHat}.

To prove \ref{ItemB2}, note that $\widehat\Delta$ can be written as
\begin{equation}
\label{Factorization}
\widehat\Delta(z) = z^4 \int_0^1 t (t-1)^2 e^{-z t} \diff t,
\end{equation}
as one immediately verifies integrating by parts. Let $s$ be a root of $\Delta$ with $s \neq s_0$. Then $z = \tau(s - s_0)$ is a root of $\widehat\Delta$ with $z \neq 0$, and thus
\begin{align}
\int_0^1 t (t - 1)^2 e^{-\sigma t} \cos(\omega t) \diff t & = 0, \label{SystIntegralCos} \displaybreak[0] \\
\int_0^1 t (t - 1)^2 e^{-\sigma t} \sin(\omega t) \diff t & = 0, \label{SystIntegralSin}
\end{align}
where $\sigma = \Real z$ and $\omega = \Imag z$. Since $\widehat\Delta$ is a quasipolynomial with real coefficients, one may assume with no loss of generality that $\omega \geq 0$.

We have $\omega > 0$, since, if $\omega = 0$, then \eqref{SystIntegralCos} is in contradiction with the fact that $t \mapsto t (t-1)^2 e^{-\sigma t}$ is strictly positive in $(0, 1)$. We also have $\omega > \pi$, since, if $\omega \leq \pi$, then $t \mapsto t (t-1)^2 e^{-\sigma t} \sin(\omega t)$ is strictly positive in $(0, 1)$. Assume, to obtain a contradiction, that $\sigma \geq 0$. Then, by Lemma~\ref{OnlyLemma}, we obtain that $\pi < \omega < 2 \pi$. By performing a suitable linear combination of \eqref{SystIntegralCos} and \eqref{SystIntegralSin}, we have
\[
\int_0^1 t (t - 1)^2 e^{-\sigma t} \sin\left(\omega \left(t - \tfrac{1}{2}\right)\right) \diff t = 0.
\]
Decomposing this integral in the intervals $(0, \frac{1}{2})$ and $(\frac{1}{2}, 1)$ and performing the changes of variables $\bar t = \frac{1}{2} - t$ and $\hat t = t - \frac{1}{2}$ in these intervals, respectively, we obtain
\[
\int_0^{\frac{1}{2}} \left(\tfrac{1}{4} - t^2\right) \sin(\omega t) \left(\sinh(\sigma t) + 2 t \cosh(\sigma t) \right) \diff t = 0,
\]
where $\sinh$ and $\cosh$ denote hyperbolic sine and cosine, respectively. Since $\pi < \omega < 2 \pi$ and $\sigma \geq 0$, the function $t \mapsto \left(\tfrac{1}{4} - t^2\right) \sin(\omega t) \left(\sinh(\sigma t) + 2 t \cosh(\sigma t) \right)$ is strictly positive in $\left(0, \frac{1}{2}\right)$, contradicting the above equality and thus concluding the proof of \ref{ItemB2}.
\end{pf*}

\section{Main result}
\label{SecMain}

The main result of the paper is the following theorem providing conditions for the maximal multiplicity of a root of $\Delta$ and establishing the dominance of this root. Recall that the maximal multiplicity of any root of $\Delta$ is $2n$.

\begin{thm}
\label{MainTheoN}
Consider the quasipolynomial $\Delta$ given by \eqref{Delta} and let $s_0 \in \mathbb R$.
\begin{enumerate}
\item\label{ItemA} The number $s_0$ is a root of multiplicity $2n$ of $\Delta$ if and only if, for every $k \in \{0, \dotsc, n-1\}$,
\begin{equation}
\label{Coeffs}
\left\{
\begin{aligned}
a_k & = \binom{n}{k} (-s_0)^{n-k} \\
& \hphantom{=} {} + (-1)^{n-k} n! \sum_{j=k}^{n-1} \binom{j}{k} \binom{2n-j-1}{n-1} \frac{s_0^{j-k}}{j!\tau^{n-j}}, \\
\alpha_k & = (-1)^{n-1} e^{s_0 \tau} \sum_{j=k}^{n-1} \frac{(-1)^{j-k} (2n-j-1)!}{k! (j-k)! (n-j-1)!} \frac{s_0^{j-k}}{\tau^{n-j}}.
\end{aligned}
\right.
\end{equation}
\item\label{ItemB} If \eqref{Coeffs} is satisfied, then $s_0$ is a strictly dominant root of $\Delta$.
\end{enumerate}
\end{thm}

By considering the first equation in \eqref{Coeffs} with $k = n-1$, one also obtains a simple expression for $s_0$ in terms of $\tau$ and $a_{n-1}$.

\begin{cor}
Let $s_0 \in \mathbb R$, $\Delta$ be the quasipolynomial given by \eqref{Delta}, and assume that the coefficients of $\Delta$ are given by \eqref{Coeffs}. Then
\[
s_0 = -\frac{a_{n-1}}{n} - \frac{n}{\tau}.
\]
\end{cor}

The complete proof of Theorem~\ref{MainTheoN} is provided in an extended version of this paper (see \cite{MazantiMultiplicity}) and follows the same lines of the proof of Theorem~\ref{MainTheo2}. After a suitable change of variables in order to reduce to the case $s_0 = 0$ and $\tau = 1$, explicit computations of $\Delta^{(k)}(z)$ for $k \in \{0, \dotsc, 2n-1\}$ allow computing the coefficients $a_0, \dotsc, a_{n-1}, \alpha_0, \dotsc, \alpha_{n-1}$ ensuring maximal multiplicity. Dominance of $s_0$ can be established by adapting the arguments from Section~\ref{SecN2}. Further insight on the sensitivity of the design can be found in \cite{Michiels2017Explicit}.

\section{Illustrative example}
\label{SecExpl}

Consider the case $n = 3$ and $\tau = 2.5$. We choose $s_0 = -0.5$ as a root of multiplicity $6$ for $\Delta$. Then, by Theorem~\ref{MainTheoN}, the coefficients $a_0, a_1, a_2, \alpha_0, \alpha_1, \alpha_2$ are given by
\begin{align*}
a_0 & = -1.735, & a_1 & = 2.91, & a_2 & = -2.1, \displaybreak[0] \\
\alpha_0 & \approx 1.736219, & \alpha_1 & \approx 1.443984, & \alpha_2 & \approx 0.3438058.
\end{align*}

\begin{figure*}[ht]
\centering
\begin{tabular}{@{} c @{} c @{}}
\resizebox{0.4494\textwidth}{!}{\input{Figures/zeros.pgf}} & \resizebox{0.4494\textwidth}{!}{\input{Figures/sols.pgf}} \tabularnewline
(a) & (b) \tabularnewline
\end{tabular}
\caption{(a) Roots of $\Delta$ for $n = 3$, $\tau = 2.5$, and $s_0 = -0.5$, within the region $\{s \in \mathbb C \suchthat -5 \leq \Real s \leq 1,\; -30 \leq \Imag s \leq 30\}$. (b) Some solutions of \eqref{MainSystTime} for $n = 3$, $\tau = 2.5$, and $s_0 = -0.5$.}
\label{FigExample}
\end{figure*}

Figure~\ref{FigExample}(a) shows the roots of $\Delta$ within the region $\{s \in \mathbb C \suchthat -5 \leq \Real s \leq 1,\; -30 \leq \Imag s \leq 30\}$, computed numerically using Python \texttt{cxroots} package, which implements numerical methods described in \cite{Kravanja2000Computing}. We observe in Figure~\ref{FigExample}(a) the presence of the multiple root $s_0 = -0.5$ as well as other roots dominated by $s_0$. Figure~\ref{FigExample}(b) illustrates the stability of \eqref{MainSystTime} with the above choice of coefficients by representing four particular solutions, obtained with the initial conditions $y_{0, i}$ for $i \in \{1, 2, 3, 4\}$ given by
\begin{align*}
y_{0, 1}(t) & = 1, & y_{0, 2}(t) & = -t, \displaybreak[0] \\
y_{0, 3}(t) & = -\frac{t^2}{4}, & y_{0, 4}(t) & = -\frac{1}{6 \omega^2} \sin(\omega t),
\end{align*}
with $\omega = 2 \pi$.

\begin{ack}
This work is supported by a public grant overseen by the French National Research Agency (ANR) as part of the ``Investissement d'Avenir'' program, through the iCODE project funded by the IDEX Paris-Saclay, ANR-11-IDEX-0003-02. The authors also acknowledge the support of Institut Polytechnique des Sciences Avanc\'ees (IPSA).
\end{ack}

\bibliography{root}

\end{document}